\documentclass[11pt,reqno]{amsart}
\usepackage{amsmath, amssymb, amsthm}
\usepackage{url}
\usepackage[breaklinks]{hyperref}
\renewcommand{\(}{\left\(}
\renewcommand{\)}{\right\)}
\renewcommand{\[}{\left\[}
\renewcommand{\]}{\right\]}

\numberwithin{equation}{section}
 \theoremstyle{plain}
\newtheorem{theorem}{Theorem}[section]

\newtheorem{remark}[]{Remark}

\newtheorem{corollary}[theorem]{Corollary}

   \makeatletter
\def\proof{\@ifnextchar[{\@oproof}{\@nproof}}
\def\@oproof[#1][#2]{\trivlist\item[\hskip\labelsep\textit{#2 Proof of\
#1.}~]\ignorespaces}
\def\@nproof{\trivlist\item[\hskip\labelsep\textit{Proof.}~]\ignorespaces}
\makeatother

\begin{document}
\title[Bressoud-Subbarao type weighted partition identities]{Bressoud-Subbarao type weighted partition identities for a generalized divisor function}

\thanks{$2020$ \textit{Mathematics Subject Classification.} Primary 11P84, 05A17; Secondary 11P81.}

\author{Archit Agarwal}
\address{Archit Agarwal, Department of Mathematics, Indian Institute of Technology Indore, Simrol, Indore 453552, Madhya Pradesh, India.}
\email{phd2001241002@iiti.ac.in, archit.agrw@gmail.com}

\author{Subhash Chand Bhoria}
\address{Subhash Chand Bhoria, Pt. Chiranji Lal Sharma Government PG College, Karnal, Urban Estate, Sector-14, Haryana 132001, India.}
\email{scbhoria89@gmail.com}

\author{Pramod Eyyunni}
\address{Pramod Eyyunni,  Department of Mathematics, Indian Institute of Technology Indore, Simrol, Indore 453552, Madhya Pradesh, India.}
\email{pramodeyy@gmail.com}

\author{Bibekananda Maji}
\address{Bibekananda Maji, Department of Mathematics, Indian Institute of Technology Indore, Simrol, Indore 453552, Madhya Pradesh, India.}
\email{bibekanandamaji@iiti.ac.in, bibek10iitb@gmail.com}

\begin{abstract}
In 1984, Bressoud and Subbarao obtained an interesting weighted partition identity for a generalized divisor function,  by means of  combinatorial arguments.  Recently,  the last three named authors found an analytic proof of the aforementioned identity of Bressoud and Subbarao starting from a $q$-series identity of Ramanujan. In the present paper,  we revisit the combinatorial arguments of Bressoud and Subbarao, and derive a more general weighted partition identity.  Furthermore,  with the help of a fractional differential operator,  we establish a few more Bressoud-Subbarao type weighted partition identities beginning from an identity of Andrews, Garvan and Liang.  
We also found a one-variable generalization of an identity of Uchimura related to Bell polynomials.  
\end{abstract}

\maketitle

\section{Introduction}

Ramanujan \cite[pp.~354--355]{ramanujanoriginalnotebook2}, \cite[pp. 302--303]{ramanujantifr}, noted down five interesting $q$-series identities at the end of his second notebook.  These identities have been recently systematically studied in \cite{BEM21}, \cite{DM}.  Here we mention one of these five $q$-series identities,  namely,  
for $ |q|<1,  c \in \mathbb{C}, 1- c q^n \neq 0, \ n \geq 1$,  
\begin{equation}\label{entry4}
\sum_{n=1}^{\infty} \frac{(-1)^{n-1} c^n q^{\frac{n(n+1)}{2} } }{(1-q^n) (c q)_n  } = \sum_{n=1}^{\infty} \frac{c^n q^n }{1-q^n}.
\end{equation}
This identity was rediscovered by Uchimura \cite[Equation (3)]{uchimura81} and Garvan \cite{garvan1},  whereas the special case $c=1$ of \eqref{entry4} goes back to Kluyver \cite{kluyver}, 
\begin{equation}\label{Kluyver}
\sum_{n=1}^{\infty} \frac{(-1)^{n-1} q^{\frac{n(n+1)}{2}}}{(1-q^n) ( q)_n  } = \sum_{n=1}^{\infty} \frac{ q^n }{1-q^n}.
\end{equation}
Uchimura  \cite{uchimura81,  uchimura87} established a new expression for the above identity with the help of a sequence of polynomials $U_1(x)=x, U_n(x)= n x^n + (1-x^n) U_{n-1}(x)$.  These polynomials have connections to the analysis of the data structure called ``heap".  Mainly,  he \cite[Theorem 2]{uchimura81}  found that 
\begin{equation}\label{Uchimura}
\sum_{n=1}^\infty n q^n (q^{n+1})_\infty =\sum_{n=1}^{\infty} \frac{(-1)^{n-1} q^{\frac{n(n+1)}{2} } }{(1-q^n) ( q)_n  } = \sum_{n=1}^{\infty} \frac{ q^n }{1-q^n}.
\end{equation}
In 1984,  Bressoud and Subbarao \cite{bresub} extracted a beautiful partition theoretic interpretation of the above identity.  Before presenting their identity,  let us first write down a few notations that will be essential throughout the article.
\begin{itemize}
\item $\pi$: an integer partition,

\item $p(n)$:= the number of integer partitions of $n$,

\item $p^{(t)}(n)$:= the number of partitions of $n$ into exactly $t$ distinct part sizes,

\item $ s(\pi):=$ the smallest part of $\pi$,

\item $\ell(\pi):=$ the largest part of $\pi$,

\item $\#(\pi):=$ the number of parts of $\pi$,

%\item $\mathrm{rank}(\pi)= \ell(\pi)- \#(\pi)$,

% \item $L(\pi):=$ total number of appearances of the largest part of $\pi$,
  
\item $\nu_{d}(\pi):=$ the number of parts of $\pi$ without multiplicity,

%\item $\nu(j):=$ the number of times the integer $j$ occurs in a partition,

%\item $\overline{p}(n):=$ the number of overpartitions of $n$,

\item $\mathcal{P}(n):=$ collection of all integer partitions of $n$,

\item $\mathcal{D}(n):=$ collection of all partitions of $n$ into distinct parts,

%\item $\mathcal{P}_{o}(n):=$ collection of all overpartitions of $n$,

%\item $\mathcal{P}^{*}(n):=$ partitions into consecutive integers with smallest part $1$.
\end{itemize}
Let $d(n)$ be the number of positive divisors of $n$.  Then comparing the coefficient of $q^n$ in \eqref{Uchimura},  one has
\begin{equation}\label{BS}
\sum_{ \pi \in \mathcal{D}(n)  } (-1)^{ \# (\pi)-1 } s(\pi)=d(n).  
\end{equation}
This identity was rediscovered by Fokking, Fokking and Wang \cite{ffw95} in 1995. 
Over the years,  the above identity motivated mathematicians to find different weighted partition identities for divisor functions.  Interested readers can see \cite{andrews08,   bresub,  dilcher,  DM, guozeng2015,  xu}.  Using purely combinatorial explanations,  Bressoud and Subbarao \cite{bresub}
found an elegant weighted partition identity for a generalized divisor function $\sigma_{z}(n)=\sum_{d|n} d^z$.  Mainly, they derived that,  for any \emph{non-negative integer} $z$, 
\begin{equation}\label{BS-general}
\sum_{ \pi \in \mathcal{D}(n)  } (-1)^{ \# (\pi)-1 } \sum_{j=1}^{s(\pi)} 
( \ell(\pi) - s(\pi) + j ) ^z = \sigma_z(n),
\end{equation}
whereas they did not point out a generating function identity for \eqref{BS-general}. 
Recently, Bhoria, Eyyunni and Maji \cite[Remark 4]{BEM21} found a generating function identity for  \eqref{BS-general} inspired from the identity \eqref{entry4} of Ramanujan and also derived a one-variable generalization of  \eqref{BS-general}.  They \cite[Theorem 2.4]{BEM21} proved that,  for any {\it integer} $z$,  complex number $c$,  $n \in \mathbb{N}$, 
\begin{equation}\label{One_Var_BS-general}
\sum_{ \pi \in \mathcal{D}(n)  } (-1)^{ \# (\pi)-1 } \sum_{j=1}^{s(\pi)} (\ell(\pi) - s(\pi)  + j ) ^z c^{\ell(\pi) - s(\pi) +j}=  \sigma_{z,c}(n),  
\end{equation}
where 
\begin{align}\label{generalized divisor function}
 \sigma_{z,c}(n)=  \sum_{d | n} d^z c^d.
 \end{align}
To derive \eqref{One_Var_BS-general},  they utilized a differential and an integral operator successively on the weighted partition implication of \eqref{entry4}.  One of the major goals of the current paper is to show that the above identity \eqref{One_Var_BS-general} is in fact true for any {\it complex number} $z$.   
%Unfortunately,  we feel that this goal cannot achieved by the method developed by Bhoria,  Eyyunni and Maji. 
 In the present paper,  we shall give detailed explanations of the combinatorial argument developed by Bressoud and Subbarao \cite{bresub} to show that \eqref{One_Var_BS-general} is indeed true for any complex number $z$, see Theorem \ref{Generalization-BS}.
 
Apart from obtaining the leftmost expression in \eqref{Uchimura}, Uchimura also generalized this identity as a whole. He defined, for each non-negative integer $m$, 
%\begin{equation*}
%U_k(q) := \sum_{n=1}^{\infty} n^k q^n (q^{n+1})_{\infty} 
%\end{equation*}
%for each positive integer $k$. He also let
\begin{equation}\label{Uch original}
M_m := M_m(q) = \sum_{n=1}^{\infty} n^m q^n (q^{n+1})_{\infty}, \quad \text{and} \quad K_{m+1} := K_{m+1}(q)= \sum_{n=1}^{\infty} \sigma_m(n) q^n.
\end{equation} 
He proved the following theorem.
\begin{theorem}\label{uch gen}
We have the following properties: \\
%\begin{itemize} 
$(1)$ $ \displaystyle \exp \left( \sum_{m=1}^{\infty} K_m t^m / m! \right) =  1+ \sum_{m=1}^{\infty} M_m t^m / m!.$ \\
$(2)$ Let $Y_m$ be the Bell polynomial defined by
\begin{equation*}
Y_m \left( u_1, u_2, \dots, u_m\right) = \sum_{\Pi(m)} \frac{n!}{k_1 ! \dots k_m !} \left( \frac{u_1}{1!} \right)^{k_1} \dots \left( \frac{u_m}{m!} \right)^{k_m},
\end{equation*}
where $\Pi(m)$ denotes a partition of $m$ with
\begin{equation*}
k_1 + 2 k_2 + \cdots + mk_m = m.
\end{equation*}
%\end{itemize}
Then for any $m \geq 1$, $M_m = Y_m (K_1, \dots, K_m)$.
\end{theorem}
Uchimura did not study the partition-theoretic implications of the above theorem, that is, he did not state the corresponding generalization of \eqref{BS}. Dilcher \cite[Corollary $1$]{dilcher} recorded a few combinatorial identities arising out of Theorem \ref{uch gen}. For $m \leq 4$, he connected the coefficients $C_m(n)$ of $q^n$ in $M_m$ with the divisor functions as follows:
\begin{align}
 C_1(n) &= d(n), \quad C_2(n) = \sigma_1(n) + \sum_{j=1}^{n-1} d(j) d(n-j), \label{dil int k=2} \\
C_3(n) &= \sigma_2(n) + 3\sum_{j=1}^{n-1} d(j) \sigma_1 (n-j) + \sum_{\substack{j+k+\ell = n \\ j, k, \ell \geq 1}} d(j)d(k)d(\ell), \nonumber \\
 C_4(n) &= \sigma_3(n) + 3\sum_{j=1}^{n-1} \sigma_1(j) \sigma_1(n-j) + 4\sum_{j=1}^{n-1} d(j) \sigma_2 (n-j) + 6 \sum_{\substack{j+k+\ell = n \\ j, k, \ell \geq 1}} d(j) d(k) \sigma_1 (\ell) \nonumber \\    & + \sum_{\substack{j_1 + \cdots + j_4 = n \\ j_1, \dots, j_4 \geq 1 }} d(j_1)\cdots d(j_4). \nonumber
\end{align} 
More importantly, the coefficient $C_m(n)$ has an interesting partition theoretic interpretation, namely,
\begin{align}
C_m(n)= \sum_{\pi \in \mathcal{D}(n)} (-1)^{\#(\pi) - 1} s(\pi)^m c^{s(\pi)}.
\end{align}
Dilcher \cite[Corollary $2$]{dilcher} derived two different generalization of \eqref{Uchimura}.  The first one is the following identity involving the binomial coefficient $\displaystyle\binom {n}{k}$. 
\begin{theorem}\label{binomial dil}
Let $k$ be a positive integer. Then
\begin{equation*}
\sum_{n=k}^{\infty} \binom {n}{k} q^n (q^{n+1})_{\infty} = q^{-{k\choose 2}  } \sum_{n=1}^\infty \frac{(-1)^{n-1} q^{n+k\choose 2}}{(1-q^n)^k (q)_n} =  \sum_{j_1 = 1}^{\infty} \frac{q^{j_1}}{1 - q^{j_1}}  \sum_{j_2 = 1}^{j_1} \frac{q^{j_2}}{1 - q^{j_2}} \cdots  \sum_{j_k = 1}^{j_{k-1}} \frac{q^{j_k}}{1 - q^{j_k}}.
\end{equation*}
\end{theorem}
Employing the above identity and Theorem \ref{uch gen} of Uchimura, Dilcher obtained the below identity. 
\begin{theorem}\label{Dilcher2}
Let $m \in \mathbb{N}$ and $K_i$'s are defined as in \eqref{Uch original}.  Then
\begin{align}
\sum_{n=1}^\infty \frac{(-1)^{n-1} q^{n+1\choose 2}}{(1-q^n)^m (q)_n} = N_m (K_1,  K_2,  \cdots,   K_m),
\end{align}
where $N_m(x_1,  x_2,  \cdots,  x_m)$ is a polynomial in $m$ variables with rational coefficients.  
\end{theorem}
Recently,  Gupta and Kumar \cite{GK21} found an interesting generalization of Theorem \ref{Dilcher2} and encountered the same generalized divisor function $\sigma_{z,c}(n)$ \eqref{generalized divisor function}.  In the same paper,  they further studied many analytic properties of $\sigma_{z,c}(n)$.  In the present paper,  we mainly focus on the weighted partition representation for $\sigma_{z,c}(n)$ and one of the other main results is a one-variable generalization of Uchimura's identity i.e.,  Theorem \ref{uch gen}.  Surprisingly,  we see the presence of the same generalized divisor function in the generalization of Theorem \ref{uch gen}.  

% (SHOULD WE ALSO MENTION ABOUT UCHIMURA IDENTITY GENERALIZATION FOR COMPLEX $m$??)   
In the next section, we mention all the main results of this paper.  
\section{Main results}

\begin{theorem}\label{Generalization-BS}
Let $c$ and $z$ be two complex numbers. For any $n \in \mathbb{N}$, we have
\begin{equation} \label{Two_Var_BS-general}
\sum_{ \pi \in \mathcal{D}(n)  } (-1)^{ \# (\pi)-1 } \sum_{j=1}^{s(\pi)} (\ell(\pi) - s(\pi)  + j ) ^z c^{\ell(\pi) - s(\pi) +j}= \sigma_{z,c}(n).
% \sum_{d | n} d^z c^d.
\end{equation}
\end{theorem}

The next theorem is a one-variable generalization of Uchimura's result in Theorem \ref{uch gen}. 

\begin{theorem}\label{gen of Uch Thm}
For a non-negative integer $m$ and a complex number $c$ with $|cq|<1$,  we define
\begin{align} \label{M_{m,c}}
M_{m,c} := M_{m,c}(q)= \sum_{n=1}^{\infty}n^mc^nq^n(q^{n+1})_{\infty}, \hspace{5mm} K_{m+1,c} :=  K_{m+1,c} (q) = \sum_{n=1}^{\infty} \sigma_{m,c}(n)q^n. \hspace{3mm}
\end{align}
%$$ 
%\textrm{where} \hspace{2mm}\sigma_{m,c}(n)=\sum_{d \mid n} d^m c^d.
%$$

Also, let $Y_m$ be the Bell polynomial defined by
\begin{align*}
Y_m (u_1,u_2,\cdots, u_m):=\sum_{\Pi(m)}\frac{m!}{k_1!\cdots k_m!}\bigg(\frac{u_1}{1!}\bigg)^{k_1}\cdots\bigg(\frac{u_m}{m!}\bigg)^{k_m},
\end{align*} 
where $\Pi(m)$ denotes a partition of $m$ with $$k_1+2k_2+\cdots + m k_m=m.$$
%\begin{itemize}
Then the exponential generating functions of $M_{m, c}$ and $K_{m,c}$ are related by:
\begin{align}
\frac{(q)_\infty}{(cq)_\infty} \exp\bigg(\sum_{m=1}^{\infty}K_{m,c}\frac{t^m}{m!}\bigg)=  \frac{(q)_\infty}{(c q)_\infty} +  \sum_{m=1}^{\infty}M_{m,c}\frac{t^m}{m!}, \quad \text{and}
\end{align}
for any $m\geq 1$, we have
\begin{align}\label{gen_Uchimura}
M_{m,c}=\frac{(q)_\infty}{(cq)_\infty}Y_m(K_{1,c},\cdots, K_{m,c}).
\end{align}
%\end{itemize}
\end{theorem}
Next, we mention a few Bressoud-Subbarao type weighted partition identities.
\begin{theorem}\label{D on AGL}
Let $k$ and $c$ be two complex numbers.  We have the following identity:
\begin{equation}\label{one var gen of FFW exp k}
\sum_{\pi \in \mathcal{D}(n)} (-1)^{\#(\pi) - 1} s(\pi)^k c^{s(\pi)} = \sum_{\pi \in \mathcal{P}(n)} \sum_{j=0}^{\nu_d(\pi)} (-1)^j \binom {\nu_d(\pi)}{j} (\ell (\pi) - j)^k c^{\ell(\pi) - j}.
\end{equation}
\end{theorem}
Note that the left hand side above arises precisely from the partition-theoretic interpretation of $M_{m, c}$ in Theorem \ref{gen of Uch Thm}. Setting $c=1$ in \eqref{one var gen of FFW exp k} naturally leads to the partition-theoretic explanation of the coefficient of $q^n$ in the series $M_m$ in \eqref{Uch original}.

\begin{corollary}\label{FFW exponent k}
Suppose that $k$ is a complex number.  Then
\begin{equation}\label{FFW gen exp k}
\sum_{\pi \in \mathcal{D}(n)} (-1)^{\#(\pi) - 1} s(\pi)^k = \sum_{\pi \in \mathcal{P}(n)} \sum_{j=0}^{\nu_d(\pi)} (-1)^j \binom {\nu_d(\pi)}{j} (\ell (\pi) - j)^k.
\end{equation}
\end{corollary}

This corollary is also a generalization of \eqref{BS}. We explain why this is so. The left hand side readily reduces to the left side of \eqref{BS} for $k=1$. The right hand side takes the form $\sum_{\pi \in \mathcal{P}(n)} \sum_{j=0}^{\nu_d(\pi)} (-1)^j \binom {\nu_d(\pi)}{j} (\ell (\pi) - j)$. This can be written as
\begin{align}
& \sum_{\pi \in \mathcal{P}(n)} \ell(\pi) \sum_{j=0}^{\nu_d(\pi)} (-1)^j \binom {\nu_d(\pi)}{j} - \sum_{\pi \in \mathcal{P}(n)} \sum_{j=0}^{\nu_d(\pi)} (-1)^j \binom {\nu_d(\pi)}{j}j \nonumber \\
& = \sum_{\pi \in \mathcal{P}(n)} \ell(\pi) (1-1)^{\nu_d(\pi)} - \sum_{\pi \in \mathcal{P}(n)} \sum_{j=0}^{\nu_d(\pi)} (-1)^j \binom {\nu_d(\pi)}{j}j \nonumber \\
&= - \sum_{\pi \in \mathcal{P}(n)} \sum_{j=0}^{\nu_d(\pi)} (-1)^j \binom {\nu_d(\pi)}{j}j.  \label{suminter}
\end{align} 
Now, for a complex number $\alpha$ and a partition $\pi$ of a positive integer, consider the identity
\begin{equation*}
(1- \alpha)^{\nu_d(\pi)} = \sum_{j=0}^{\nu_d(\pi)} (-\alpha)^j \binom {\nu_d(\pi)}{j}.
\end{equation*}
Differentiating with respect to $\alpha$, we get
\begin{equation}\label{D 1 time}
-\nu_d(\pi) (1 - \alpha)^{\nu_d(\pi) - 1} = \sum_{j=0}^{\nu_d(\pi)}(-1)^j j \alpha^{j-1} \binom {\nu_d(\pi)}{j}.
\end{equation}
We let $\alpha$ approach $1$, and see that the right side above reduces to the inner sum in \eqref{suminter}. The behaviour of the left side above depends on the value of $\nu_d(\pi)$. Note that, as $\alpha \rightarrow 1$, $-\nu_d(\pi) (1 - \alpha)^{\nu_d(\pi) - 1} $ tends to $0$ if $\nu_d(\pi) > 1$ and to $- \nu_d(\pi) = -1$ when $\nu_d(\pi) = 1$. Hence, putting this information in \eqref{suminter}, we can write
\begin{align}
\sum_{\pi \in \mathcal{P}(n)} \sum_{j=0}^{\nu_d(\pi)} (-1)^j \binom {\nu_d(\pi)}{j} (\ell (\pi) - j) &= - \sum_{\pi \in \mathcal{P}(n)} \sum_{j=0}^{\nu_d(\pi)} (-1)^j \binom {\nu_d(\pi)}{j}j \nonumber \\
&= - \sum_{\substack{\pi \in \mathcal{P}(n) \\ \nu_d(\pi) = 1}} -1 = \sum_{\substack{\pi \in \mathcal{P}(n) \\ \nu_d(\pi) = 1}} 1,  \label{suminter2}
\end{align}
the number of partitions $\pi$ of $n$ with $\nu_d(\pi) = 1$. But the partitions of $n$ with only one distinct part are of the form $m+m+ \cdots + m$, where $m$ is a positive divisor of $n$, and conversely, for each positive divisor $m$ of $n$, we do have the partition $m + \cdots +m$ with $n/m$ summands, which is a partition with one distinct part. Thus, the last sum in \eqref{suminter2} equals $d(n)$, and we have showed that \eqref{FFW gen exp k} is indeed a generalization of \eqref{BS}. 

Recall that $p^{(2)} (n)$ denotes the number of partitions of $n$ with exactly two distinct parts sizes.  
%From the identities for $C_1(n)$ and $C_2(n)$ in \eqref{dil int k=2}, and from Theorem \ref{binomial dil}.
% (SHOULD WE INCLUDE IT??) \\
From the identity \eqref{dil int k=2} of Dilcher,  for $k=2$, namely, 
\begin{equation}\label{dilcher k=2}
\sum_{\pi \in \mathcal{D}(n)} (-1)^{\#(\pi) - 1} s(\pi)^2 = \sigma(n) + \sum_{j=1}^{n-1} d(j) d(n-j),
\end{equation}
and the special case $k=2$ of \eqref{FFW gen exp k}, we obtain a representation for $p^{(2)} (n)$. 
\begin{corollary}\label{rep for exactly 2 distinct}
For each positive integer $n$, we have
\begin{equation}\label{exactly 2 distinct}
p^{(2)} (n) = \frac{1}{2} \left\{ \sum_{j=1}^{n-1}d(j)d(n-j) + d(n) - \sigma(n) \right\}.
\end{equation}
\end{corollary}
%--------------------------------------------------------------------------------------------------------------
%\begin{example}
%Let $n=6$. The partitions of $6$ with exactly two distinct parts are $5+1, \ 4+2, \ 4+1+1, \ 3+1+1+1, \ 2+2+1+1, \ 2+1+1+1+1$, and thus $p_d^{(2)} (6) = 6$. The right hand side of \eqref{exactly 2 distinct} for $n=6$ reads as 
%\begin{align*}
%\frac{1}{2} \left\{ \sum_{j=1}^{5} d(j)d(6-j) + d(6) - \sigma(6) \right\} &= \frac{1}{2} \left\{d(1)d(5) + d(2)d(4) + d(3)d(3) + d(4)d(2) + d(5)d(1) + d(6) - \sigma(6) \right\} \\ 
%&= \frac{1}{2} \left\{2 + 6 + 4 + 6 + 2 + 4 - 12 \right\} = 6.
%\end{align*}  
%\end{example}

%(SHOULD THE NEXT THM BE RETAINED??)
Again, we mention a Bressoud-Subbarao type weighted partition identity for the generalized divisor function $\sigma_{k,c}(n)$. 
\begin{theorem}\label{ACEM} Let $n \in \mathbb{N}$ and $k,  c \in \mathbb{C}$.  Then we get 
\begin{align}\label{ACEM equation}
\sum_{\pi\in \mathcal{D}(n)}(-1)^{\#(\pi)-1}\sum_{j=1}^{s(\pi)}j^k c^j= \sum_{\substack{\pi \in \mathcal{P}(n) \\ \nu_d(\pi) \geq 2}} \sum_{j=0}^{\nu_d(\pi)-1}(-1)^j {\nu_d(\pi)-1 \choose j}(\ell(\pi)-j)^k c^{\ell(\pi)-j} + \sigma_{k,c}(n).
\end{align}
\end{theorem}
\begin{remark}
Corresponding to $k=0$ and $c=1$, the above identity reduces to \eqref{BS}. 
Note that $d(n)= p^{(1)}(n)$, where $p^{(1)}(n)$ denotes the number of partitions of $n$ into exactly $1$ part size.   Hence,  identity  \eqref{BS} can be rewritten as 
\begin{align} \label{different form_BS}
p^{(1)}(n)=  \sum_{ \pi \in \mathcal{D}(n)  } (-1)^{ \# (\pi)-1 } s(\pi).  
\end{align}
\end{remark}

%
%\begin{corollary}
%(CAN WE GET FROM SIMPLER VERSION?) In the case of $k=1$ and $c=1$ we have a new representation for the sum of divisors,
%\begin{align} 
%\sum_{\pi\in \mathcal{D}(n)}(-1)^{\#(\pi)-1}{s(\pi)+1 \choose 2}=\sum_{d|n}d\quad + \quad \sum_{\substack{\pi \in \mathcal{P}(n)  \\ \nu_d(\pi) =2}}1
%\end{align}
%\end{corollary}

Letting $k=1$ and $c=1$ in Theorem \ref{ACEM}, we obtain an interesting analog of \eqref{different form_BS}. 
\begin{corollary}\label{Analogue of BS}
Let $n$ be a positive integer. Then we have
%(TRY COMBINATORIAL PROOF, CAN WE GET FROM SIMPLER VERSION)
\begin{align}\label{Analogue of BS equation}
 p^{(2)}(n) = \sum_{\pi \in \mathcal{D}(n)} (-1)^{ \# (\pi) }s(\pi) \left( \ell(\pi) - s(\pi) \right). 
\end{align} 
% where $\rho(\pi) =\ell(\pi) - s(\pi)$,
% call it  interval length of the partition in which all the parts lie.
\end{corollary}

%Next, we state another beautiful  weighted partition identity, namely, (NO NEED TO STATE THIS AS A SEPARATE RESULT??) 
%\begin{theorem}
%\begin{align}
%\sum_{\pi\in \mathcal{D}(n)}(-1)^{\#(\pi)-1}(\ell(\pi)+1)s(\pi)=2\sum_{d|n}d\quad + \quad \sum_{\substack{\pi \in \mathcal{P}(n)  \\ \nu_d(\pi) =2}}1
%\end{align}
%\end{theorem}
In the next section,  we demonstrate proofs of all main results. 
 
\section{Proof of the main results}
\begin{proof}[Theorem \ref{Generalization-BS}][]
The main idea of the proof of this theorem is due to Bressoud and Subbarao.  Here we explain the details while simultaneously extending their result to complex number $z$.  
%We adapt the idea Bressoud and Subbarao used in proving \eqref{BS-general}, to prove this more general result.  
We first define, for each positive integer $N$, the set $\mathcal{C}(N)$ of partitions $\pi$ into distinct parts satisfying the following inequalities:
\begin{align*}
    \ell(\pi) \geq N>\ell(\pi)-s(\pi).
\end{align*}
%For example,  one can check that $\mathcal{C}(1)$ contains the partitions $\{ (1),  (2),  (3),  \cdots \}$. 
 Let $\pi \in \mathcal{C}(N)$ be any partition into distinct parts.  The definition of $\mathcal{C}(N)$ implies that the only possibilities for $N$ are $\ell(\pi)-s(\pi)+j$, with $1\leq j\leq s(\pi)$. Therefore for any partition $\pi$ into distinct parts, there are exactly $s(\pi)$ many integers $N$ such that $\pi \in \mathcal{C}(N)$.  With the help of this fact the left hand side of \eqref{Two_Var_BS-general} can be written as 
\begin{align}\label{left side}
\sum_{ \pi \in \mathcal{D}(n)  } (-1)^{ \# (\pi)-1 } \sum_{j=1}^{s(\pi)} (\ell(\pi) - s(\pi)  + j ) ^z c^{\ell(\pi) - s(\pi) +j} &  = \sum_{\pi \in \mathcal{D}(n)} (-1)^{\#(\pi)-1}\sum_{\substack{N=1\\
\pi\in \mathcal{C}(N)}}^{n}N^z c^N  \nonumber \\ &=\sum_{N=1}^{n} N^z c^N \sum_{\substack{\pi\in \mathcal{D}(n)\\ \pi\in \mathcal{C}(N)}}(-1)^{\#(\pi)-1} .  
\end{align}
Now our goal is to show that 
\begin{align}
 \sum_{\substack{\pi\in \mathcal{D}(n)\\ \pi\in \mathcal{C}(N)}}(-1)^{\#(\pi)-1} = \begin{cases} 
      1, & {\rm if} \, N \mid n, \\
      0, & {\rm if}\, N \nmid n. 
   \end{cases} \label{to show}
\end{align}
To prove \eqref{to show}, we shall try to pair the distinct parts partitions of $n$ in $\mathcal{C}(N)$ that have opposite parity in their respective number of parts. The contribution of such pairs will be a $+1$ and a $-1$ to the sum in the left side of \eqref{to show}, since the weight is $(-1)^{\#(\pi)-1}$. Thus, the sum of the contributions of each such pair will vanish.  When $N\mid n$, the only distinct parts partition that will remain is $n$ itself.  And when $N\nmid n$,  all the distinct parts partitions of $n$ will be paired.  

Now we shall explain the algorithm to create such pairs. We consider two different cases.
\textbf{Case 1$\colon$} Let $\pi \in \mathcal{D}(n) \cap \mathcal{C}(N)$. If $\pi$ contains a part which is divisible by $N$ and if $\pi$ has at least one other part, (note that the other part cannot be a multiple of $N$ since $ \ell(\pi)-s(\pi) < N$) then remove the part which is divisible by $N$ and add $N$ to the smallest remaining part.  
Continue this way to create new partitions by adding $N$ to the smallest part in the previous partition until we get a distinct parts partition of $n$.

A natural question arises here as to what the guarantee is that this algorithm gives a new distinct parts partition of $n$ that also lies in $\mathcal{C}(N)$. In the next paragraph,  we shall explain by taking a general distinct parts partition of $n$.  

Let $\pi \in \mathcal{D}(n) \cap \mathcal{C}(N)$ with parts $a_1<a_2<\cdots<a_k$ and suppose $N$ divides $a_i$ for some $1\leq i\leq k$,  say $a_i=j\cdot N$ for some $j \geq 1$.   First,  we  remove $a_i$ and add $N$ to the smallest remaining part.  We continue adding $N$ to the smallest part in the previous partition to create a new partition until we again have a partition of $n$.  Since $a_i = j \cdot N$,  we have to do this process $j$ times.  Let $\pi_j'$ denote the new partition of $n$.  
%{\bf Is it okay if we denote it simply by $\pi_j$ ?}

Here we look at various cases depending on the value of $j$. If $j<k-1$,  then $\pi_j'$ equals
\begin{align*}
\begin{cases}
a_{j+1}+\cdots+a_{i-1}+a_{i+1}+\cdots+a_{k}+(a_{1}+N)+\cdots +(a_{j}+N), & {\rm if}\,  1\leq j<i-1, \\
a_{j+2}+\cdots +a_{k}+(a_{1}+N)+\cdots+(a_{i-1}+N)+(a_{i+1}+N)+\cdots +(a_{j+1}+N),  & {\rm if}\, i-1\leq j<k-1.
\end{cases} \end{align*}
One can easily check that the new partition $\pi_j'$ is a distinct parts partition of $n$ and lies in $\mathcal{C}(N)$.  Now, if $j=k-1$,  then we simply add $N$ to each part and get a new distinct parts partition of $n$,  namely,  
$\pi_j'= (a_1 + N) + \cdots + (a_{i-1}+N)+(a_{i+1}+N)+\cdots +(a_{k}+N)$, which also belongs to $\mathcal{C}(N)$.  Note that $\pi$ and $\pi_j'$ have opposite parity in their number of parts.  Again, if $j > k-1$,  then with the help of division algorithm and utilizing previous cases, one can construct $\pi_j' \in \mathcal{D}(n) \cap \mathcal{C}(N)$ such that $\pi$ and $\pi_j'$ will have opposite parity in their number of parts. 

%  { \bf In this process,  we are able to pair $\{\pi,  \pi_j'\}$ such that they have opposite parity in number of parts. }
 
 \textbf{Case 2$\colon$} Let $\pi \in \mathcal{D}(n) \cap \mathcal{C}(N)$ with parts $a_1<a_2<\cdots<a_k$ such that $N \nmid a_i$ for all $1 \leq i \leq k$.  In this case,  we must reverse the procedure,  that is, we subtract $N$ from the largest part in the previous partition until we reach a unique partition $\pi'$ for which \begin{align}
\ell(\pi')-N<\textrm{ the total amount subtracted} <s(\pi')+N. \label{stopping condition}
\end{align} 
Finally, this total amount subtracted is then inserted as a new part.
 The above condition \eqref{stopping condition} may look artificial,  but soon we will explain why this condition comes in naturally.  

First,  we subtract $N$ from the largest part $a_k$. Then we have $a_k - N < a_1$ since $ N > \ell(\pi) - s(\pi) = a_k -a_1$.  
Let $\pi_1:= (a_k - N) + a_1 + \cdots + a_{k-1}$.  Note that $\pi_1$ is not a partition of $n$,  so we need to add $N$ to $\pi_1$ to get a distinct parts partition of $n$. This brings up three possibilities.

\textbf{Sub case 1}: Firstly, suppose that $N < a_k-N$.  Let us define $\pi_1':= N + (a_k-N) + a_1 + \cdots + a_{k-1}$. See that $\pi_1' \in \mathcal{D}(n)$. Moreover, $\pi_1' \in \mathcal{C}(N)$ if and only if $ \ell(\pi_1')= a_{k-1} \geq N > \ell(\pi_1') - s(\pi_1') =  \ell(\pi_1')- N$. In this sub case, the first inequality will be satisfied naturally since we assumed that $N < a_k - N$,  but the second inequality may or may not be true.  The second inequality will be true if  $N > \ell(\pi_1')- N (\textrm{the amount subtracted})$,  that is,  $N (\textrm{the amount subtracted}) > \ell(\pi_1')- N $.

\textbf{Sub case 2}: Next, suppose we have $a_k - N < N < a_{k-1}$. In this sub case, let us define $\pi_1':=  (a_k-N) +  \cdots +N + \cdots + a_{k-1}$,  which is a distinct parts partition of $n$.  Now $\pi_1' \in \mathcal{C}(N)$ if and only if $ \ell(\pi_1')= a_{k-1} \geq N > \ell(\pi_1') - s(\pi_1') = a_{k-1} - (a_k - N)$.  In this sub case, one can easily see that both of the inequalities are true. 

\textbf{Sub case 3}: Finally, suppose $a_{k-1} < N$ is true. Now we define $\pi_1':=  (a_k-N) + a_1 + \cdots + a_{k-1} + N$, which is a distinct parts partition of $n$.  Here $N$ is the largest part.  Now $\pi_1' \in \mathcal{C}(N)$ if and only if  $\ell(\pi_1') \geq N > \ell(\pi_1') - s(\pi_1') = N - s( \pi_1')$.  The first inequality is true since $\ell(\pi_1') = N$, but the second inequality will be true if  $s(\pi_1') + N > N (\textrm{the amount subtracted})$. 
Combining all three sub cases,  we can clearly see that $\pi_1' \in \mathcal{C}(N)$ if and only if 
\begin{align}\label{1st}
 \ell(\pi_1')- N < N \, (\textrm{the amount substracted}) < s(\pi_1') + N.
 \end{align}
Thus,  the above condition justifies why we need \eqref{stopping condition}.  
If the chain of inequalities in \eqref{1st} does not hold at this stage,  then again we subtract $N$ from the largest part $a_{k-1}$ in $\pi_1$.  And we define $\pi_2:= (a_{k-1}-N) + (a_k - N) + a_1 + \cdots + a_{k-2}$.  Now one can observe that we have subtracted $2N$ from the original partition $\pi$ of $n$,  so we have to add $2N$ to obtain a distinct parts partition of $n$.  Again, we face three sub cases.  Firstly,  if $2 N < (a_{k-1}-N)$,  then we write $\pi_2':= 2N + (a_{k-1}-N)  + \cdots + a_{k-2}$.  Secondly, if $ (a_{k-1}-N) < 2N < a_{k-2}$,  then put $\pi_2':=  (a_{k-1}-N) + \cdots+ 2N + \cdots + a_{k-2}$. 
And thirdly,  if $a_{k-2} < 2N$,  then set $\pi_2':=  (a_{k-1}-N)  + \cdots + a_{k-2}+ 2N$.
Along the same lines, we can show that $\pi_2' \in \mathcal{C}(N)$ if and only if 
\begin{align}\label{2nd}
 \ell(\pi_2')- N < 2N \, (\textrm{the amount subtracted}) < s(\pi_2') + N.
 \end{align}
If the partition $\pi_2'$ does not satisfy the above conditions, then we apply the algorithm once again.  More generally,
suppose that we have subtracted the integer $N$ $`j$' many times. We then have to add $j \cdot N$ as a new part. Note that this part may be the smallest part, the largest part or neither of them, and so we have to consider three subcases.   If $j<k$, then our new distinct parts partition of $n$, again denoted by $\pi_j'$, will be
 \begin{align*}
\begin{cases} j \cdot N+(a_{k-j+1}-N)+\cdots+a_k+\cdots+a_{k-j}, \,\, \textrm{if} \,\, j \cdot N\, \textrm{ is the smallest  part},  \\ 
(a_{k-j+1}-N)+\cdots+j \cdot N+ \cdots+a_{k-j},  \,\,  \textrm{if} \,\,  j \cdot N \, \textrm{ is neither the smallest nor the largest part}, \\ 
(a_{k-j+1}-N)+\cdots+a_k+\cdots+a_{k-j}+j \cdot N,  \,\,  \textrm{if}\,\,  j \cdot N \, \textrm{is the largest part}.  \end{cases}
\end{align*} 
In a similar vein,  one can show that  $\pi_j' \in \mathcal{C}(N)$ if and only if  
\begin{align}
\ell(\pi_j')-N< j \cdot N \, (\textrm{the total amount subtracted}) <s(\pi_j')+N.
\end{align} 
This justifies why condition \eqref{stopping condition} is necessary.  If $j \geq k$,  then we can use division algorithm and give similar arguments. Hence, summarizing Case 2, we started with a distinct parts partition $\pi$ with no part divisible by $N$ and then constructed another distinct parts partition $\pi_j'$ which has only one part divisible by $N$,  namely,  $j \cdot N$ for some suitable positive integer $j$.  

The above algorithm explains that if $N \nmid n$,  then any partition $\pi \in \mathcal{D}(n) \cap \mathcal{C}(N)$ can be paired with another partition $\pi_j' \in \mathcal{D}(n) \cap \mathcal{C}(N)$ such that they have opposite parity in their number of parts.  And if $N |n$,  the only partition of $n$ which will remain unpaired is the partition $n$ itself since the partition $n$ neither belongs to Case 1 nor to Case 2.   This completes the proof of \eqref{to show}.  Finally,  combining \eqref{left side} and \eqref{to show},  one can obtain \eqref{Two_Var_BS-general}. 
\end{proof}

\begin{remark}
In \cite[Section 5]{BEM21}, the authors showed that Theorem \ref{Generalization-BS} is valid for any integer $z$, by applying the differential operator $D[f(c)]:=c \displaystyle \frac{d}{dc}\{f(c)\}$, and the integral operator
$I[f(c)]:= \displaystyle \int_{0}^{c}\frac{f(t)}{t}dt$, on the partition theoretic interpretation of  \eqref{entry4}.   Here we have  stated that Theorem \ref{Generalization-BS} is in fact true for any complex number $z$. Moreover,  for  positive real numbers $z$ and $j$,  one can define the fractional derivative
\begin{align*}
\frac{d^z}{{dc}^z}(c^j):= \frac{\Gamma(j+1)}{\Gamma(j-z+1)} c^{j-z}. 
\end{align*}
This definition matches with the usual definition of the derivative when $z$ is any positive integer. Therefore, using this fractional derivative, we have $D^{z}(c^j)= j^z c^j$ for any $z>0$. Hence applying this fractional operator on  the partition theoretic interpretation \cite[Lemma 5.1]{BEM21} of  \eqref{entry4},  one can first show that, Theorem \ref{Generalization-BS} is true for any positive real number $z$ and then using analytic continuation, we can prove that the identity \eqref{Two_Var_BS-general} is in fact valid for any complex $z$. 
\end{remark}

\begin{proof}[Theorem \ref{gen of Uch Thm}][]
Let us define a function
\begin{equation}\label{G(x,q)}
A(x,q) :=\frac{(q)_\infty}{(xq)_\infty}.
\end{equation} 
We next invoke Euler's identity, which is a special case of the $q$-binomial theorem. For $|q|,  |t| < 1$, we have
\begin{equation}\label{Euler q-binom}
\sum_{n=0}^{\infty} \frac{t^n}{(q)_n}=\frac{1}{(t)_\infty}.
\end{equation}
Replacing $t$ by $x q$ in \eqref{Euler q-binom} and  utilizing in \eqref{G(x,q)}, we get,  for $|xq| < 1$,
\begin{equation*}
A(x,q) =(q)_\infty\sum_{n=0}^{\infty}\frac{x^nq^n}{(q)_n} = \sum_{n=0}^{\infty} x^nq^n(q^{n+1})_\infty. 
\end{equation*} 
Putting $x=ce^t,$   for $m \geq 1$,  one can see that
\begin{align*}
    \frac{\partial^m }{\partial t^m}  A(ce^t,q)\bigg|_{t=0} :=\sum_{n=1}^{\infty}n^mc^nq^n(q^{n+1})_{\infty} =M_{m,c}.  
\end{align*}
Thus,  the power series for $A(c e^t,  q)$ in the variable $t$ takes the shape as
\begin{align}\label{power series for G}
A(c e^t,  q) = A(c, q) + \sum_{m=1}^\infty M_{m,c} \frac{t^m}{m!}.  
\end{align}
Suppose we write
\begin{align}\label{log(G)}
\log(A(x,q))= \sum_{n=0}^{\infty}h_n(c) \frac{t^n}{n!}.
\end{align} From \eqref{G(x,q)}, we get
\begin{align*}
    \log(A(x,q))=\log((q)_\infty)-\log((xq)_\infty) &=\log((q)_\infty)-\sum_{n=1}^{\infty}\log(1-xq^n)\\
    &=\log((q)_\infty)+\sum_{n=1}^{\infty}\sum_{m=1}^{\infty}\frac{(xq^n)^m}{m}\\
    &=\log((q)_\infty)+\sum_{m=1}^{\infty}\frac{x^m}{m}\frac{q^m}{1-q^m}.
\end{align*} Putting $x=ce^t,$ we find that
\begin{align*}
    \log(A(ce^t,q))&=\log((q)_\infty)+\sum_{m=1}^{\infty} \frac{q^m}{1-q^m}\frac{c^me^{tm}}{m}\\
    &=\log((q)_\infty)+\sum_{m=1}^{\infty}\frac{q^m}{1-q^m}\frac{c^m}{m}\sum_{n=0}^{\infty}\frac{(tm)^n}{n!}\\
    &=\log((q)_\infty)+\sum_{n=0}^{\infty}\sum_{m=1}^{\infty}\frac{c^mm^{n-1}q^m}{1-q^m}\frac{t^n}{n!}\\
    &=\log((q)_\infty)+\sum_{m=1}^{\infty}\frac{c^m}{m}\frac{q^m}{1-q^m}+\sum_{n=1}^{\infty}\bigg(\sum_{m=1}^{\infty}\frac{c^mm^{n-1}q^m}{1-q^m}\bigg)\frac{t^n}{n!}\\
    &=\log((q)_\infty)-\log((cq)_\infty)+\sum_{n=1}^{\infty}\bigg(\sum_{m=1}^{\infty}\frac{c^mm^{n-1}q^m}{1-q^m}\bigg)\frac{t^n}{n!}\\
    &=\log\bigg(\frac{(q)_\infty}{(cq)_\infty}\bigg)+\sum_{n=1}^{\infty}\bigg(\sum_{m=1}^{\infty}\frac{c^mm^{n-1}q^m}{1-q^m}\bigg)\frac{t^n}{n!}.
\end{align*}By comparing coefficients in \eqref{log(G)}, we get
\begin{align}\label{h_0(c)}
h_0(c)=\log\bigg(\frac{(q)_\infty}{(cq)_\infty}\bigg),
\end{align}
and for any $n\geq0$, we have
\begin{align}\label{h_n(c)}
h_{n+1}(c)=\sum_{m=1}^{\infty}\frac{c^mm^nq^m}{1-q^m}=\sum_{m=1}^{\infty}\sum_{k=1}^{\infty}c^m  m^nq^{mk}=\sum_{\ell=1}^{\infty}\sum_{d\mid \ell}c^dd^nq^\ell=K_{n+1,c}.
\end{align}
Now in view of \eqref{log(G)},  \eqref{h_0(c)},  and \eqref{h_n(c)},  one can see that
\begin{align*}
    A(ce^t,q)=\frac{(q)_\infty}{(cq)_\infty} \exp\bigg(\sum_{n=1}^{\infty}K_{n,c}\frac{t^n}{n!}\bigg),
\end{align*} and finally using \eqref{power series for G},  we arrive at
\begin{align*}
    \frac{(q)_\infty}{(cq)_\infty} \exp\bigg(\sum_{n=1}^{\infty}K_{n,c}\frac{t^n}{n!}\bigg)=  \frac{(q)_\infty}{( c q)_\infty} + \sum_{n=1}^{\infty}M_{n,c}\frac{t^n}{n!}.
\end{align*}
Finally comparing the coefficients of $t^n/n!$,  for $n\geq 1$, and using the definition of the Bell polynomial, one can obtain \eqref{gen_Uchimura}.  
\end{proof}

\begin{proof}[Theorem \ref{D on AGL}][]
We start with the partition-theoretic interpretation of an identity due to Andrews, Garvan and Liang \cite{agl13}, first noted down by Dixit and Maji \cite[Equation $(2.6)$, Corollay $2.4$]{DM}
\begin{equation}\label{PTI of AGL}
\sum_{\pi \in \mathcal{D}(n)} (-1)^{\#(\pi) - 1} \left( 1 + c + \cdots + c^{s(\pi) - 1}\right) = \sum_{\pi \in \mathcal{P}(n)} c^{\ell(\pi) - \nu_d(\pi)} (c-1)^{\nu_d(\pi) - 1}.
\end{equation}
Multiplying throughout by $(c-1)$, we get an identity with which we are going to work.
\begin{equation}\label{final AGL PTI} 
\sum_{\pi \in \mathcal{D}(n)} (-1)^{\#(\pi) - 1} \left(c^{s(\pi)} - 1\right) = \sum_{\pi \in \mathcal{P}(n)} c^{\ell(\pi) - \nu_d(\pi)} (c-1)^{\nu_d(\pi)}.
\end{equation}
The left side is clearly a polynomial in $c$. The same is true of the right side as well, as $\ell(\pi) \geq \nu_d(\pi)$ for any partition $\pi$. To see why, observe that if $\ell(\pi)$ is the largest part in a partition, then the possible parts in the partition come from the set $\{ 1, 2, \dots, \ell(\pi)\}$. This means that $\nu_d(\pi)$, the number of distinct parts appearing in the partition $\pi$ is at most $\ell(\pi)$, with equality occuring if and only if each of the integers from $1$ to $\ell(\pi)$ appears in the partition. We now apply the fractional differential operator $D^k:= \left( c \frac{d}{{dc}} \right)^k$, with $k>0$, to both sides of \eqref{final AGL PTI}. The left side clearly transforms into
\begin{equation}\label{LHS after D}
\sum_{\pi \in \mathcal{D}(n)} (-1)^{\#(\pi) - 1} s(\pi)^k \ c^{s(\pi)}.
\end{equation}
Before applying $D^k$ to the right side of \eqref{final AGL PTI}, we expand it using the binomial theorem to get
\begin{equation}\label{RHS poly}
\sum_{\pi \in \mathcal{P}(n)} c^{\ell(\pi) - \nu_d(\pi)} \sum_{j=0}^{\nu_d(\pi)} \binom {\nu_d(\pi)}{j} (-1)^j c^{\nu_d(\pi) - j} = \sum_{\pi \in \mathcal{P}(n)} \sum_{j=0}^{\nu_d(\pi)} c^{\ell(\pi) - j} (-1)^j \binom {\nu_d(\pi)}{j}.
\end{equation}
Operating by $D^k$, we obtain
\begin{equation}\label{RHS after D}
\sum_{\pi \in \mathcal{P}(n)} \sum_{j=0}^{\nu_d(\pi)} (\ell(\pi) - j)^k c^{\ell(\pi) - j} (-1)^j \binom {\nu_d(\pi)}{j}.
\end{equation}
Equating the differentiated expressions in \eqref{LHS after D} and \eqref{RHS after D}, we get \eqref{one var gen of FFW exp k} for any $k>0$. Finally, making use of analytic continuation on the variable $k$, the proof of Theorem \ref{D on AGL} is over.
\end{proof}

%(AS COMMENTED IN THE MAIN RESULTS SECTION, SHOULD THE FOLLOWING COROLLARY BE RETAINED SINCE THE RESULT CAN BE DERIVED FROM DILCHER'S AND UCHIMURA'S RESULTS??) \\
%----------------------------------------------------------------------------------------------------------------- \\
\begin{proof}[Corollary \ref{rep for exactly 2 distinct}][]
The idea of this proof is to consider the case $k=2$ of \eqref{FFW gen exp k} in Corollary \ref{FFW exponent k}, and compare it with \eqref{dilcher k=2} due to Dilcher. Begin by setting $k=2$ in \eqref{FFW gen exp k} to see that
\begin{align}
\sum_{\pi \in \mathcal{D}(n)} (-1)^{\#(\pi) - 1} s(\pi)^2 &= \sum_{\pi \in \mathcal{P}(n)} \sum_{j=0}^{\nu_d(\pi)} (-1)^j \binom{\nu_d(\pi)}{j} (\ell(\pi) - j)^2 \nonumber \\
&=  \sum_{\pi \in \mathcal{P}(n)} \sum_{j=0}^{\nu_d(\pi)} (-1)^j \binom{\nu_d(\pi)}{j} \left( \ell(\pi)^2 - 2 \ell(\pi)j + j^2 \right) \nonumber \\
&= \sum_{\pi \in \mathcal{P}(n)} \ell(\pi)^2 (1-1)^{\nu_d(\pi)} -2 \sum_{\pi \in \mathcal{P}(n)} \ell(\pi) \sum_{j=0}^{\nu_d(\pi)} (-1)^j \binom{\nu_d(\pi)}{j} j \nonumber  \\
& \quad + \sum_{\pi \in \mathcal{P}(n)} \sum_{j=0}^{\nu_d(\pi)} (-1)^j \binom{\nu_d(\pi)}{j} j^2 \nonumber \\
&=  -2 \sum_{\pi \in \mathcal{P}(n)} \ell(\pi) \sum_{j=0}^{\nu_d(\pi)} (-1)^j \binom{\nu_d(\pi)}{j} j + \sum_{\pi \in \mathcal{P}(n)} \sum_{j=0}^{\nu_d(\pi)} (-1)^j \binom{\nu_d(\pi)}{j} j^2. \label{2 expressions}  
\end{align}
In the first sum in \eqref{2 expressions}, we have already 
%(REFER to the MAIN RESULTS section, if REQD label this as a COROLLARY)
 seen in \eqref{suminter} that  the inner sum equals $-1$ if $\nu_d(\pi) = 1$, and $0$ if $\nu_d(\pi) > 1$. So the first sum in \eqref{2 expressions} reduces to $2 \displaystyle \sum_{\substack{ \pi \in \mathcal{P}(n) \\ \nu_d(\pi) = 1}} \ell(\pi)$. We have previously observed that partitions of a positive integer $n$ with one distinct part correspond to divisors of $n$. This means that if $m + \cdots + m$ is a partition of $n$, corresponding to a divisor $m$ of $n$, then the largest part is also $m$. Thus,   
\begin{equation}\label{simplified 1 distinct}
2 \sum_{\substack{ \pi \in \mathcal{P}(n) \\ \nu_d(\pi) = 1}} \ell(\pi) = 2 \sum_{d \ | \ n} d. 
\end{equation}
Putting this in \eqref{2 expressions} brings us to 
\begin{equation}\label{intermediate}
\sum_{\pi \in \mathcal{D}(n)} (-1)^{\#(\pi) - 1} s(\pi)^2 = 2 \sum_{d \ | \ n} d + \sum_{\pi \in \mathcal{P}(n)} \sum_{j=0}^{\nu_d(\pi)} (-1)^j \binom{\nu_d(\pi)}{j} j^2.
\end{equation}
Recall the identity \eqref{D 1 time} and multiply it with $c$ to get 
\begin{equation*}
-\nu_d(\pi) \cdot c \cdot (1-c)^{\nu_d(\pi) - 1} = \sum_{j=0}^{\nu_d(\pi)} (-1)^j j c^j \binom{\nu_d(\pi)}{j}.
\end{equation*}
Now differentiate this with respect to $c$ to derive the identity
\begin{equation}\label{D 2 times}
-\nu_d(\pi) \left[ (1-c)^{\nu_d(\pi) - 1} - c (\nu_d(\pi) - 1) (1-c)^{\nu_d(\pi) - 2} \right] = \sum_{j=0}^{\nu_d(\pi)} (-1)^j \cdot j^2 \cdot c^{j-1} \cdot \binom{\nu_d(\pi)}{j}.  
\end{equation}
Letting $c \rightarrow 1$ in \eqref{D 2 times}, we see that 
\begin{equation}\label{final}
\sum_{j=0}^{\nu_d(\pi)} (-1)^j \cdot j^2 \cdot \binom{\nu_d(\pi)}{j} = F(\nu_d(\pi)),
\end{equation}
where $$F(\nu_d(\pi)) = \lim_{c \rightarrow 1} -\nu_d(\pi) \left[ (1-c)^{\nu_d(\pi) - 1} - c (\nu_d(\pi) - 1) (1-c)^{\nu_d(\pi) - 2} \right].$$
If $\nu_d(\pi) > 2$, $F(\nu_d(\pi))$ clearly approaches $0$. When $\nu_d(\pi) = 2$, the first term in the parenthetical sum in $F(\nu_d(\pi))$ vanishes and the second term nears the value $-1$, so that in effect, $F(2) = 2$. Finally, when $\nu_d(\pi) = 1$, the second term in $F(\nu_d(\pi))$ goes to $0$ and the first term to $1$, thereby giving us the value of $F(1)$ as $-1$. In summary,
$$
F(\nu_d(\pi)) =
\begin{cases}
-1, \quad \text{if} \ \nu_d(\pi) = 1,  \\
2, \quad \text{if} \ \nu_d(\pi) = 2, \\
0, \quad \text{if} \ \nu_d(\pi) > 2.
\end{cases}
$$
Using this, via \eqref{final}, in \eqref{intermediate}, we arrive at
\begin{equation}\label{final eqn for pd2}
\sum_{\pi \in \mathcal{D}(n)} (-1)^{\#(\pi) - 1} s(\pi)^2 = 2 \sum_{d \ | \ n} d - \sum_{\substack{\pi \in \mathcal{P}(n) \\ \nu_d(\pi) = 1}} 1 + \sum_{\substack{\pi \in \mathcal{P}(n) \\ \nu_d(\pi) = 2}} 2 = 2 \sigma(n) - d(n) +  2 p_d^{(2)}(n).
\end{equation}  
Comparing \eqref{final eqn for pd2} with \eqref{dilcher k=2}, we see that
\begin{equation}
\sigma(n) + \sum_{j=1}^{n-1} d(j) d(n-j) = 2 \sigma(n) - d(n) +  2 p_d^{(2)}(n),
\end{equation}
which on rearrangement yields \eqref{exactly 2 distinct} and the proof is complete.
\end{proof}

\begin{proof}[Theorem \ref{ACEM}][]
We start with \eqref{PTI of AGL} and multiplying throughout by $c$, we get
\begin{equation}\label{c-AGL PTI}
\sum_{\pi \in \mathcal{D}(n)} (-1)^{\#(\pi) - 1} \sum_{j=1}^{s(\pi)} c^j = \sum_{\pi \in \mathcal{P}(n)} c^{\ell(\pi) - \nu_d(\pi) + 1} (c-1)^{\nu_d(\pi)-1}.
\end{equation} Now apply the fractional differential operator $D^k := \left( c \frac{d}{{dc}} \right)^k$, with $k>0$,  to both sides of \eqref{c-AGL PTI}. The left side clearly transforms into \begin{equation}\label{LHS of Thm 2.6}
\sum_{\pi \in \mathcal{D}(n)} (-1)^{\#(\pi) - 1} \sum_{j=1}^{s(\pi)} j^k c^j.
\end{equation} Before applying the fractional differential operator to the right hand side of \eqref{c-AGL PTI}, we use binomial theorem to expand
\begin{equation}\label{RHS c-AGL PTI}
\sum_{\pi \in \mathcal{P}(n)} c^{\ell(\pi) - \nu_d(\pi) + 1} (c-1)^{\nu_d(\pi)-1} = \sum_{\pi \in \mathcal{P}(n)} \sum_{j=0}^{\nu_d(\pi)-1} (-1)^j \binom {\nu_d(\pi)-1}{j}c^{\ell(\pi) - j}.
\end{equation} 
Now applying the fractional differential operator $D^k$, we obtain
\begin{align}
\sum_{\pi \in \mathcal{P}(n)} \sum_{j=0}^{\nu_d(\pi)-1} (-1)^j \binom {\nu_d(\pi)-1}{j} ( \ell(\pi) - j)^k c^{\ell(\pi) - j} & = \sum_{\substack{\pi \in \mathcal{P}(n) \\ \nu_d(\pi) \geq 2}}  \sum_{j=0}^{\nu_d(\pi)-1}(-1)^j {\nu_d(\pi)-1 \choose j}(\ell(\pi)-j)^k c^{\ell(\pi)-j}\nonumber \\
&+ \sum_{\substack{\pi \in \mathcal{P}(n) \\ \nu_d(\pi)=1}} \ell(\pi)^k. \label{RHS of Thm 2.6}
\end{align} 
We have previously observed that partitions of a positive integer $n$ with one distinct part correspond to divisors of $n$. This means that if $m + \cdots + m$ is a partition of $n$, corresponding to a divisor $m$ of $n$, then the largest part is also $m$. Thus,
\begin{equation}
\sum_{\substack{\pi \in \mathcal{P}(n) \\ \nu_d(\pi)=1}} \ell(\pi)^k = \sum_{d \in n} d^k
\end{equation} Putting this in \eqref{RHS of Thm 2.6}, we get
\begin{equation}\label{final RHS of Thm 2.6}
 \sum_{\substack{\pi \in \mathcal{P}(n) \\ \nu_d(\pi) \geq 2}} \sum_{j=0}^{\nu_d(\pi)-1}(-1)^j {\nu_d(\pi)-1 \choose j}(\ell(\pi)-j)^k c^{\ell(\pi)-j} + \sum_{d|n}d^k c^d
\end{equation}
Equating the differentiated expressions in \eqref{LHS of Thm 2.6} and \eqref{final RHS of Thm 2.6}, we can obtain \eqref{ACEM equation} for $k>0$. Finally, using analytic continuation on $k$, one can see that
 Theorem \ref{ACEM} is valid for any complex $k$.
\end{proof}

\begin{proof}[Corollary \ref{Analogue of BS}][]
Putting $k=1$ and $c=1$ in Theorem \ref{ACEM},  the identity becomes
\begin{align}\label{k=c=1}
\sum_{\pi\in \mathcal{D}(n)}(-1)^{\#(\pi)-1}\sum_{j=1}^{s(\pi)}j = \sum_{\substack{\pi \in \mathcal{P}(n) \\ \nu_d(\pi) \geq 2}} \sum_{j=0}^{\nu_d(\pi)-1}(-1)^j {\nu_d(\pi)-1 \choose j}(\ell(\pi)-j) + \sum_{d|n}d . 
\end{align} 
Letting $z=c=1$ in Theorem \ref{Generalization-BS},   we have
\begin{align}\label{z=c=1}
\sum_{ \pi \in \mathcal{D}(n)  } (-1)^{ \# (\pi)-1 } \sum_{j=1}^{s(\pi)} (\ell(\pi) - s(\pi)  + j ) = \sum_{d | n} d.
\end{align}
From \eqref{k=c=1} and \eqref{z=c=1},  one can see that
\begin{align}
\sum_{\pi \in \mathcal{D}(n)} (-1)^{ \# (\pi) }s(\pi) \left( \ell(\pi) - s(\pi) \right) =  \sum_{\substack{\pi \in \mathcal{P}(n) \\ \nu_d(\pi) \geq 2}} \sum_{j=0}^{\nu_d(\pi)-1}(-1)^j {\nu_d(\pi)-1 \choose j}(\ell(\pi)-j). 
\end{align}
Now we will divide the right hand sum into two parts corresponding to $\nu_d(\pi)=2$ and $\nu_d(\pi) \geq 3$.  Thus,  we have
\begin{align}\label{two parts}
\sum_{\substack{\pi \in \mathcal{P}(n) \\ \nu_d(\pi) = 2}} \sum_{j=0}^{1}(-1)^j {1 \choose j}(\ell(\pi)-j) + \sum_{\substack{\pi \in \mathcal{P}(n) \\ \nu_d(\pi) \geq 3}} \sum_{j=0}^{\nu_d(\pi)-1}(-1)^j {\nu_d(\pi)-1 \choose j}(\ell(\pi)-j).  
\end{align}
One can easily see that the first sum reduces to $p^{(2)}(n)$.  To show that the second sum vanishes, 
%$$
%\sum_{\substack{\pi \in \mathcal{P}(n) \\ \nu_d(\pi) = 2}} 1 = p^{(2)}(n).  
%$$
%\begin{equation}\label{3.28}
%\sum_{\pi\in \mathcal{D}(n)}(-1)^{\#(\pi)} s(\pi)(\ell(\pi)-s(\pi))=\sum_{\substack{\pi \in \mathcal{P}(n) \\ \nu_d(\pi) = 2}} 1 + \sum_{\substack{\pi \in \mathcal{P}(n) \\ \nu_d(\pi) \geq 3}} \sum_{j=0}^{\nu_d(\pi)-1}(-1)^j {\nu_d(\pi)-1 \choose j}(\ell(\pi)-j)
%\end{equation} 
we mention the following  expressions:
\begin{equation}\label{binomial}
\begin{aligned}
(1-x)^m&= \sum_{j=0}^{m} (-1)^j {m \choose j} x^j,  \\
-m(1-x)^{m-1}&= \sum_{j=0}^{m} (-1)^j {m \choose j}j x^{j-1}. 
\end{aligned}
\end{equation} Now letting $x=1$ and $m=\nu_d(\pi)-1$ in \eqref{binomial} and using them,  one can see that the second sum in \eqref{two parts} vanishes.  
\end{proof}

\section{Concluding Remarks}

In this paper,  we established an interesting weighted partition identity \eqref{Two_Var_BS-general} for a generalized divisor function $\sigma_{z, c}(n):= \sum_{d|n} d^z c^d$.  Mainly,  we adopt the combinatorial proof of Bressoud and Subbarao.  We also point out that Theorem \ref{Generalization-BS} can also be derived from the partition theoretic interpretation of \eqref{entry4} by using a fractional differentiation operator.  Motivated from this observation,  we applied the aforementioned operator on the partition theoretic interpretation of an identity of Andrews,  Garvan and Liang,  and obtained a few Bressoud-Subbarao type weighted partition identities.  
% which is not a common technique in the literature. 

We also generalize Uchimura's identity i.e.,  Theorem \ref{uch gen} by introducing a complex parameter $c$ and surprisingly, the generating functions $K_n$ for the divisor functions $\sigma_{z, 1}(n)$ in Theorem \ref{uch gen} are replaced by the generating functions for the generalized divisor functions $\sigma_{z, c}(n)$,  for $z \in \mathbb{N}$.
%appearing in \eqref{One_Var_BS-general}. 
%(see Theorem \ref{gen of Uch Thm} in the next section) 
This raises a natural question.  Since the generalization of the identity of Bressoud and Subbarao in \eqref{One_Var_BS-general} holds true for complex numbers $z$ as well,  it would be interesting to see if such a generalization exists for Theorem \ref{uch gen}. (i.e., an analogue of Theorem \ref{gen of Uch Thm},  for complex numbers $m$).  

In \eqref{different form_BS},  we mentioned an alternate form of Bressoud-Subbarao's identity \eqref{BS}. 
We found an interesting analogue of \eqref{different form_BS} i.e.,  Corollary \ref{Analogue of BS},  which is a weighted partition identity for $p^{(2)}(n)$.  It would be highly desirable to find a Bressoud-Subbarao type combinatorial proof for Corollary \ref{Analogue of BS}.

\section{Acknowledgements}

The second author wants to thank the Department of Mathematics,
Pt. Chiranji Lal Sharma Government College, Karnal for a conducive research environment.
The third author is a SERB National Post Doctoral Fellow (NPDF) supported by the
fellowship PDF/2021/001090 and would like to thank SERB for the same. The last author
wishes to thank SERB for the Start-Up Research Grant SRG/2020/000144.  

\section{Data availability statement}
The authors declare that the manuscript has no associated data. 

\section{Statement of Conflict of Interest}

On behalf of all authors,  the corresponding author states that there is no conflict of interest.

\end{document}